\documentclass[11pt,reqno]{amsart}
\usepackage{amsmath,amssymb,amsfonts}
\usepackage{eucal,enumerate}

\newcommand\vstrut[2]{\rule[#1]{0ex}{#2}}

\numberwithin{equation}{section}

\newtheorem{lem}{Lemma}[section]
\newtheorem{cor}[lem]{Corollary}

\newtheorem{thm}[lem]{Theorem}

\theoremstyle{remark}
\newtheorem{exam}[lem]{Example}
\newtheorem{prob}[lem]{Problem}
\newtheorem{conj}[lem]{Conjecture}

\newcommand\ts{\textstyle}

\renewcommand{\phi}{\varphi}                 
\renewcommand{\epsilon}{\varepsilon}
\newcommand\eset{\emptyset}                  
\renewcommand\emptyset{\varnothing}
\newcommand\subeq{\subseteq}

\newcommand\inv{^{-1}}
\newcommand\Span[1]{\langle #1 \rangle} 
\newcommand\aff{\operatorname{aff}}
\newcommand\vol{\operatorname{vol}}

\newcommand\cH{\mathcal{H}}
\newcommand\cL{\mathcal{L}}
\newcommand\cP{\mathcal{P}}
\newcommand\cS{\mathcal{S}}

\newcommand\bbR{\mathbb{R}}
\newcommand\bbZ{\mathbb{Z}}

\newcommand\bj{\mathbf{1}}

\newcommand\setm{\setminus}
\newcommand\0{\hat 0}

\newcommand\G{\Gamma}

\newcommand\PH{{P,\cH}}	
\newcommand\PoH{{P^\circ,\cH}}	
\newcommand\Orth{O}

\renewcommand\qedsymbol{\ensuremath{\blacksquare}}

\renewcommand\aa{{\mathrm{a}}}
\newcommand\cc{{\mathrm{c}}}

\newcommand\Tehrhyp{Theorem 4.2}	
\newcommand\Tquasi{Theorem 4.1}		

\begin{document}

\begin{center}
\Large
An Enumerative Geometry for Magic and Magilatin Labellings
\normalsize
\vskip20pt

Matthias Beck \\ Department of Mathematics\\ San Francisco State University\\ 1600 Holloway Avenue\\ San Francisco, CA 94132, U.S.A.\\
E-mail: {\tt beck@math.sfsu.edu}\\[10pt]

and\\[10pt]

Thomas Zaslavsky\\
Department of Mathematical Sciences \\
State University of New York at Binghamton\\
Binghamton, NY 13902-6000, U.S.A.\\
E-mail: {\tt zaslav@math.binghamton.edu}\\[20pt]

{Version of \today.}\\[20pt]

\end{center}

\Small{\sc Abstract.}
 A \emph{magic labelling} of a set system is a labelling of its points by 
distinct positive integers so that every set of the system has the same 
sum, the \emph{magic sum}.  Examples are magic squares (the sets are the 
rows, columns, and diagonals) and semimagic squares (the same, but without 
the diagonals).  A \emph{magilatin labelling} is like a magic labelling 
but the values need be distinct only within each set.  We show that the 
number of $n\times n$ magic or magilatin labellings is a quasipolynomial 
function of the magic sum, and also of an upper bound on the entries in 
the square.  Our results differ from previous ones because we require that 
the entries in the square all be different from each other, and because we 
derive our results not by \emph{ad hoc} reasoning but from a general 
theory of counting lattice points in rational inside-out polytopes.  We 
also generalize from set systems to rational linear forms.  
\\



\normalsize

\section{It's all kinds of magic}  \label{intro}

We offer a theory for counting magic squares and their innumerable relatives: semimagic and pandiagonal magic squares, magic cubes and hypercubes, magic graphs, and magical oddities like circles, spheres, and stars \cite{Andrews,Ball}.

A \emph{magic square} is an $n \times n$ array of distinct positive integers whose sum along any \emph{line} (row, column, or main diagonal) is the same number, the \emph{magic sum}.
Magic squares date back to China in the first millenium B.C.E.\ \cite{CammannO}, came in the first millenium C.E.\ to the Islamic world and India \cite{CammannI}, and passed to Europe in the later Middle Ages \cite{CammannI} and to sub-Saharan Africa not much after \cite{AfCounts}.  
The contents of a magic square have varied with time and writer; usually they have been the first $n^2$ consecutive positive integers (\emph{standard} squares), but often any arithmetic sequence and sometimes fairly arbitrary numbers.  The fixed ideas are that they are integers, positive (or rarely, nonnegative), and distinct.  Even the mathematical treatises \cite{Andrews,Ball,Benson} take positivity and distinctness so much for granted as never to mention them.  

In the last century mathematicians took an interest in results about the number of squares with a fixed magic sum, but with simplifications.  Thus diagonal sums were usually omitted and, most significantly, the fundamental requirement of distinctness was almost invariably neglected.  
Our work, however, follows tradition by adhering to the distinctness requirement.  We shall call a square \emph{magic} or \emph{strongly magic} if its entries are distinct and \emph{weakly magic} if they need not be distinct; and strongly or weakly \emph{semimagic} if the diagonals are ignored.  

The literature on exact formulas examines the functions $W(t)$ and $W_0(t)$ that count positive and nonnegative weak semimagic or magic squares with magic sum $t$.  These functions are amenable to analysis in terms of Ehrhart theory \cite{Ehr}.  (See for example \cite{Esquares,LHDE} for semimagic squares; for magic squares see \cite{Bmagic}.)  
One treats a square with upper bound or magic sum $t$ as an integer vector $x \in t\cdot[0,1]^{n^2}$, confined to the subspace $t s$, where $t=1,2,3,\ldots$ and 
$$
s := \{x \in \bbR^{n^2} : \text{ all line sums equal } 1\},
$$
 the \emph{magic subspace}.  (Exactly which subspace this is depends on 
whether we treat squares that are semimagic, magic, pandiagonal magic 
[where we include among the line sums the wrapped diagonals], or of 
another type.)  Thus a square $x$ is an integer lattice point in the 
$t$-dilate $tP$ of the \emph{magic polytope} $P := [0,1]^{n^2} \cap s$; 
moreover, $x \in t P^\circ$, the relative interior of $P$, if and only if 
the square is positive.  

Ehrhart's fundamental results on integer-point enumeration in polytopes \cite{Ehr} give much insight.  Ehrhart theory implies that $W$ and $W_0$ are 
quasipolynomials in $t$.  (A \emph{quasipolynomial} is a function 
$Q(t)=\sum_0^d c_it^i$ with coefficients $c_i$ that are periodic functions 
of $t$, so that $Q$ is a polynomial on each residue class modulo some 
integer, called the \emph{period}; these polynomials are the 
\emph{constituents} of $Q$.)  The quasipolynomials are polynomials (that 
is, the periods are 1) in the semimagic case, because the matrix that 
defines $s$ is totally unimodular so the vertices of $P$ are all integral.  
In the magic case this is unfortunately not so and the period is not easy 
to calculate.

Still there were no exact (theoretical) formulas for strong squares (not even in the comprehensive tome \cite{SGPNRmagic}), with the exception of Stanley's \cite[Exercise 4.10]{EC1}.  With the theory of \emph{inside-out polytopes} \cite{IOP} we can attack this and related counting problems in a systematic way obtaining a general result about magic counting functions and an interpretation of reciprocity that leads to a new kind of question about permutations.
In inside-out theory we supplement the polytope $P=[0,1]^{n^2} \cap s$ with the pair-equality hyperplane arrangement 
$$
\cH := \cH[K_{n^2}]^s = \{ h_{ij} \cap s : i < j \leq {n^2}\},
$$
 where $h_{ij}$ is the hyperplane $x_i=x_j$, $\cH[\G] := \{ h_{ij}: ij \in E \}$ is the hyperplane arrangement of the graph $\G$ with edge set $E$, and $K_d$ denotes the complete graph on $d$ nodes.  
The number of $n\times n$ squares corresponding to $s$ with magic sum $t$ is the number of integer points in $t \left( P^\circ \setminus \bigcup \cH \right)$.  
This is a quasipolynomial in $t$ by the general theory of inside-out 
polytopes.  Then inside-out reciprocity \cite{IOP} gives the enumeration 
of weak nonnegative squares with multiplicity; this is reminiscent of Stanley's 
theorem on acylic orientations \cite{AOG}.

Another famous kind of square is latin squares and their relatives.  Here 
each line has $n$ different numbers.  In a \emph{latin square} these $n$ 
numbers are the same in every line and are normally taken to be the first 
$n$ positive integers.  In any latin square in this broad meaning, every 
line has the same sum.  Suppose we add this property to the definition of 
a latin square but we loosen the restriction on the entries, so that the 
square is filled with positive integers having equal row 
and column sums.  We call such squares \emph{magilatin}.  Then a magilatin 
square is a point in $\bbZ^{n^2}$; the only difference between a semimagic 
and a magilatin strong square is that we assume fewer inequations between 
the entries; while in a semimagic square each entry must differ from every 
other, in a magilatin square it must differ only from those that are 
collinear with it, a \emph{line} being a row 
or column. As with magic and semimagic squares, inside-out polytope theory 
yields theorems about the number of magilatin squares as a function either 
of the magic sum or of the largest allowed value of an entry in the 
square.

There is a parallel generalization of latin rectangles.  A \emph{latin 
rectangle} is an $m \times n$ rectangular array filled by $n$ symbols, 
none repeated in a row or column.  The asymptotic numbers of latin squares 
and rectangles of given dimensions have been the subject of many studies.  
Our geometric counting method leads in a different direction that, as far 
as we know, has not been studied.  Define a \emph{magilatin rectangle} to 
be a point in $\bbZ^{mn}$ but not in $\bigcup \cH[ K_m \times K_n ]$.  We 
discuss the number of magilatin squares or rectangles of fixed dimensions 
as we vary the maximum permitted value.

The magic and latin properties generalize far beyond squares and rectangles.  Semimagic and pandiagonal magic squares suggest a general picture: that of a \emph{covering clutter}, consisting of a finite set $X$ of points together with a family $\cL$ of subsets, called \emph{lines} for no particular reason, of which none contains any other and none is empty, and whose union is $X$.  We want to assign positive integers to $X$ so that all line sums are equal to a single number. Such a labelling is called a weakly or strongly magic or latin labelling of the covering clutter, depending on the particular requirements. 
 There are many interesting examples that we cannot treat individually here.  $(X,\cL)$ may be a finite affine or projective geometry, the ``lines'' being the subspaces of any fixed dimension; more generally it may be a block design.  It may be an $n \times n \times \cdots \times n$ hypercubical array.  It may be a 
\emph{$k$-net}, where the lines fall into $k$ parallel classes (with 
$k\geq2$) so that each point belongs to a unique line in each parallel 
class.  (A semimagic square is a 2-net and a pandiagonal square is a kind 
of 4-net.)  All these examples have lines of equal size, a property that 
has advantages but is not necessary for the theory to apply.

These ideas generalize still further.  Take rational linear forms 
$f_1,f_2,\ldots,f_m$.  A magic labelling of $[d]:= \{ 1, 2, \dots, d \}$ 
with respect to $f_1,f_2,\ldots,f_m$ is an integer point $x \in \bbR^d$ 
such that all the values $f_i(x)$ are equal to the same number.
The analog here of a covering clutter in which all lines have the same size is a system of forms for which all values $f_i(\bj)$ (the \emph{weights}; $\bj$ is the vector of all ones) are equal; such systems have nice properties.

We treat two distinct interesting approaches to enumeration.  Traditionally, magic and semimagic squares have been counted with the magic sum as the 
parameter (due to its geometrical interpretation we call this 
\emph{affine} counting); but another tack is to take as parameter the 
maximum allowed value of a label (which we call \emph{cubical} counting).  
These same two systems apply to latinity.  In our treatment we develop 
both counting systems equally.

The reader may wonder how practical our counting method is.  We believe it 
is relatively feasible.  In \cite{SLS} we apply it to solve in utter detail six 
problems of $3\times3$ squares: magic, semimagic, and magilatin (all strong), counted both cubically and affinely.

\section{Inside-out polytopes take the stage}\label{prelim}

The theory of inside-out polytopes \cite{IOP} is designed to count those points of the integral lattice $\bbZ^d$ that are contained in a rational convex polytope $P$ but not in a rational affine hyperplane arrangement $\cH$, that is, where each hyperplane is spanned by the rational points it contains.
 We call $(\PH)$ a \emph{rational inside-out polytope}, \emph{closed} if 
$P$ is closed.  We shall always assume $P$ is closed, except when we 
specifically state otherwise.  The affine span $\aff P$ may be a proper 
subspace of $\bbR^d$.

A \emph{region} of $\cH$ is a connected component of $\bbR^d \setm \bigcup\cH$. 
A \emph{closed region} is the closure of a region.  A \emph{region} of $(\PH)$ is the nonempty intersection of a region of $\cH$ with $P$.  
A \emph{vertex} of $(\PH)$ is a vertex of any such region.  
Note that a closed region of $(\PH)$ is the closure of an open region of 
$(\PH)$ and therefore meets the relative interior $P^\circ$.  The 
\emph{denominator} of $(\PH)$ is the smallest positive integer $t$ for 
which $t\inv\bbZ^d$ contains every vertex of $(\PH)$.

The fundamental counting functions associated with $(\PH)$ are the \emph{closed} \emph{Ehrhart quasipolynomial}, 
\begin{align*}
E_\PH(t) &:= \sum_{x\in t\inv\bbZ^d}m_\PH(x),
\intertext{where $P$ is closed and the \emph{multiplicity} $m_\PH(x)$ of $x\in \bbR^d$ with respect to $\cH$ is defined through}
m_\PH(x) &:= \begin{cases}
\text{the number of closed regions of $(\PH)$ that contain $x$}, 
  & \text{if } x \in P, \\
0, & \text{if } x \notin P,
\end{cases}
\intertext{and the \emph{open Ehrhart quasipolynomial}, }
E^\circ_\PoH(t) &:= \# \left( t\inv\bbZ^d \cap \left[ P^\circ \setminus {\ts\bigcup}\cH \right] \right) .
\end{align*}

We denote by $\vol P$ the volume of $P$ normalized with respect to $\bbZ^d 
\cap \aff P$, that is, we take the volume of a fundamental domain of 
$\bbZ^d \cap \aff P$ to be 1.  When $P$ is full dimensional this is the 
ordinary volume.

\begin{thm}[{\cite[\Tquasi]{IOP}}] \label{T:quasi}
If $(\PH)$ is a closed, full-dimensional, rational inside-out polytope in $\bbR^d$, then $E_\PH(t)$ and $E^\circ_\PoH(t)$ are quasipolynomials in $t$ that satisfy the reciprocity law $E^\circ_\PoH(t) = (-1)^d E_\PH(-t).$, with period equal to a divisor of the denominator of $(\PH)$, with leading term $(\vol P) t^d$, and with constant term $E_\PH(0)$ equal to the number of regions of $(\PH)$.
\end{thm}

In particular, if $(\PH)$ is integral then $E_\PH$ and $E^\circ_\PoH$ are polynomials. 

The \emph{M\"obius function} of a finite partially ordered set $S$ is the function $\mu: S \times S \to \bbZ$ defined recursively by 
\begin{equation*}
\mu(r,s) := \begin{cases}
	0			&\text{if } r \not\leq s, \\
	1			&\text{if } r = s, \\
	-\sum_{r \leq u < s} \mu(r,u)	&\text{if } r < s.
	\end{cases}
\end{equation*}
 Sources are, \emph{inter alia}, \cite{FCT} and \cite[Section 3.7]{EC1}. 

The \emph{intersection semilattice} of $\cH$ is defined as 
$$
\cL(\cH) := \big\{ {\ts\bigcap} \cS : \cS \subeq \cH \text{ and } {\ts\bigcap} \cS \neq \eset \big\} ,
$$ 
ordered by reverse inclusion (so the whole space $\bbR^d$, the intersection of no hyperplanes, is the bottom element $\0$).  Its members are the \emph{flats} of $\cH$.  
The \emph{intersection poset} of $(\PoH)$ is defined as 
$$
\cL(\PoH) := \big\{ P^\circ \cap {\ts\bigcap} \cS : \cS \subeq \cH \big\} \setminus \big\{ \eset \big\} ,
$$ 
ordered by reverse inclusion.  $\cL(\PoH)$ is a ranked poset and every interval is a geometric lattice.  
It equals $\cL(\cH)$ if $\bigcap\cH$ meets $P^\circ$.

The arrangement \emph{induced} by $\cH$ on $s \subseteq \bbR^d$ is
$$
\cH^s := \{ h\cap s : h\in\cH, \ h\not\supseteq s \}.
$$

For the second theorem 
we need the notion of transversality.  $\cH$ is \emph{transverse} to $P$ if every flat $u\in \cL(\cH)$ that intersects $P$ also intersects $P^\circ$, and $P$ does not lie in any of the hyperplanes of $\cH$. Let 
$$
E_P(t):= \# \left( tP \cap \bbZ^d \right),
$$ 
the standard Ehrhart counting function (without any hyperplanes present).

\begin{thm}[{\cite[\Tehrhyp]{IOP}}] \label{T:ehrhyp} 
If $P$ and $\cH$ are as in Theorem \ref{T:quasi}, then
\begin{equation}\label{E:oehrhyp}
E^\circ_\PoH(t) = \sum_{u\in \cL(\PoH)} \mu(\0,u) E_{P^\circ\cap u}(t) ,
\end{equation}
and if $\cH$ is transverse to $P$,
\begin{equation} \label{E:ehrhyp}
E_\PH(t) = \sum_{u\in \cL(\PoH)} |\mu(\0,u)| E_{P\cap u}(t) .
\end{equation}
\end{thm}

Often the polytope is not full-dimensional.  Suppose that $s$ is any affine subspace.  Its \emph{period} $p(s)$ is the smallest positive integer $p$ for which $p\inv\bbZ^d$ meets $s$.

\begin{cor}[{\cite[Corollary 4.3]{IOP}}] \label{C:affine}
Let $P$ be a rational convex polytope and $\cH$ a hyperplane arrangement in $s := \aff P$.  
Then $E_\PH(t)$ and $E^\circ_\PoH(t)$ are quasipolynomials in $t$ that satisfy the reciprocity law $E^\circ_\PoH(t) = (-1)^{\dim s}E_\PH(-t)$. 
Their period is a multiple of $p(s)$ and a divisor of the denominator of $(\PH)$.  
If $t \equiv 0 \mod {p(s)}$, the leading term of $E_\PH(t)$ is $(\vol_{p(s)\inv\bbZ^d} P) t^{\dim s}$ and its constant term is the number of regions of $(\PH)$; but if $t \not\equiv 0 \mod {p(s)}$, then $E_\PH(t) = E^\circ_\PoH(t) = 0$.
\hfill$\qedsymbol$
\end{cor}

\section{Magic squares, magic labellings of covering clutters, and equal line sums}  \label{magic}

\subsection{Sorts of magic} \label{intromagic}

We pointed out in the introduction that the difference between weak and 
strong magic (or semimagic) squares lies in the fact that for the latter 
we require the entries to be distinct.  We will therefore spend the 
beginning of this section studying the general setting of integer points 
in polytopes with distinct entries.

We have a convex polytope $P \subseteq \bbR^d$, spanning an affine 
subspace $s$.  To ensure distinctness of the coordinates of a vector we avoid the hyperplanes of $\cH := \cH[ K_d ]^s$, the 
complete-graph arrangement $\cH[ K_d ]$ intersected with $s$.  
Transversality of $\cH$ and $P$ means that, first of all, $s$ is not a subspace of any 
hyperplane $x_j = x_k$ and, secondly, any flat of $\cH[ K_d ]$ that meets 
$P$ also meets $P^\circ$.  In many of the interesting special cases the 
latter condition is automatic.

Suppose $x$ is a point in $\bbR^d$ whose entries are all distinct.  There 
is a unique permutation $\tau$ of $[d]$ such that $x_{\tau 1} < x_{\tau 2} 
< \cdots < x_{\tau d}$.  We say $x$ \emph{realizes} $\tau$.  We call a 
permutation $\sigma$ \emph{realizable} in a subset $A \subseteq s$ if 
there is a vector $x \in A$ that realizes it.  We are interested in 
realizability in $P$, but realizability in $s$ is simpler.  Fortunately, a 
permutation that is realizable in $s$ is also realizable in $P$ if $P$ 
contains a positive multiple of $\bj := (1,1,\ldots,1)$, since every 
closed region of $\cH$ contains $\Span{\bj}$.  This is the case when every 
form has equal positive weight.

If $x \in \bbR^d$ satisfies $x_{\sigma 1} \leq x_{\sigma 2} \leq \cdots \leq x_{\sigma d}$, we say $x$ and $\sigma$ are \emph{compatible}.

\begin{thm}\label{T:magic} 
Suppose $P \subseteq \bbR^d$ is a closed, rational convex polytope transverse to $\cH[ K_d ]$ and $s := \aff P$.  
The number $E^\circ_{P^\circ, \cH[K_d]^s} (t)$ of integer points in $t P^\circ$ with distinct entries is a quasipolynomial in positive integers $t$ with leading term $(\vol P) t^{\dim s}$.  Furthermore, 
$(-1)^{\dim s} E^\circ_{P^\circ, \cH[K_d]^s} (-t) = E_{P, \cH[K_d]^s} (t):=$ the number of pairs $(x,\sigma)$ consisting of an integer point $x \in tP$ and a compatible $P$-realizable permutation $\sigma$ of $[d]$. 
 The constant term of $E_{P, \cH[K_d]^s} (t)$ equals the number of permutations of $[d]$ that are realizable in $P$.
\end{thm} 

\begin{proof} 
The first statement is a direct consequence of Corollary \ref{C:affine} along with the observation that, by transversality, a region that intersects $P$ must also intersect $P^\circ$. 
For the second, the multiplicity $m(x)$ of $x$ equals the number of closed regions of $\cH[K_d]^s$ that contain $x$.  
The regions of $\cH[K_d]^s$ correspond to certain regions of $\cH[K_d]$, which correspond to permutations of $[d]$.  
Clearly, a closed region contains $x$ if and only if its permutation is compatible with $x$.  
Thus, $m(x)$ is the number of permutations that are both realizable in $P$ and compatible with $x$.  Now appeal to Corollary \ref{C:affine}. 
\end{proof} 

\begin{prob} \label{Pr:perdenom}
The period and denominator present a puzzle.  The denominator of the inside-out polytope $(P,\cH[K_d])$ is obviously a multiple of the denominator of the standard polytope $P$.  The first question is when the hyperplane arrangement $\cH[ K_d ]$ changes this latter denominator, and in what way.  
As for the period, if in particular $P$ has integral vertices then $E_P$ is a polynomial.  What conditions on $P$ ensure that $E_{P, \cH[ K_d ]}$ is also a polynomial? 
That it need not be is illustrated by the simple example of the line segment from $(0,1)$ to $(1,0)$ in $\bbR^2$ and the hyperplane $x_1 = x_2$. 
\end{prob}

Our desire to develop the ideas behind magic squares and graphs suggests two approaches to choosing $P$ and $s$.  The subspace, $s$, represents the existence of a magic sum.  The polytope, $P$, represents the constraints on the entries in the value vector.  The magic sum constraints may be pure equalities:
\begin{enumerate} 
\item[(i)] Set all linear forms equal to each other.  (Homogeneous equations.)
\end{enumerate} 
Or, they may be set all equal to a controlled constant:
\begin{enumerate} 
\item[(ii)] Set all linear forms equal to $t$.  (Affine equations.)
\end{enumerate} 
(Sometimes one wants additional equations; see the discussion of centrally symmetric squares in Examples \ref{X:homosym} and \ref{X:affsym}.)  Similarly, the constraints on the components of $x$ may be two-sided bounds:
\begin{enumerate} 
\item[(I)] All variables $x_i$ satisfy $0 \leq x_i \leq t$.  (Cubical constraints.)
\end{enumerate} 
Or, the constraints may be merely nonnegativity of the variables:
\begin{enumerate} 
\item[(II)] All variables $x_i \geq 0$.  (Nonnegativity.)
\end{enumerate} 
We think the natural combinations are (i) with (I) and (ii) with (II) and that is how we develop the theory.

\subsection{Cubical magic} \label{homomagic}

The cubical approach to magic squares counts them by the largest allowed value of the entry in a cell; if $t$ is the parameter, the squares counted are those with entries $0 < y_{ij} < t$.  Similarly, magic labellings of a bidirected graph are counted by the upper bound $t-1$ on the edge labels.

In the general situation the magic subspace $s$ is defined by homogeneous, rational linear equations 
\begin{equation} \label{E:homoeqns}
f_1(x) = f_2(x) = \cdots = f_m(x) = 0.
\end{equation}
These are obviously equivalent to (i), if the forms in \eqref{E:homoeqns} are the differences of the forms of (i).
The magic polytope is $$P := [0,1]^d \cap s,$$ whence the name 
``cubical''.  
The hyperplane arrangement is $\cH := \cH[K_d]^s$. 
In the example of magic squares, $P$ is the set of all square matrices with real entries in $[0,1]$ such that all row, column, and diagonal sums are equal.  For magic graphs $P$ is the set of edge labellings by numbers in $[0,1]$ such that all node sums are equal.

One has to be sure that $P$ spans $s$.  

\begin{lem} \label{L:homospan}
If $P$ is not contained within a coordinate hyperplane, it affinely spans the magic subspace $s$.
\end{lem}

\begin{proof}
By the hypothesis $P$ contains points $x^i$ (not necessarily distinct) with $0 < (x^i)_i \leq 1$.  The barycenter of these points lies in $(0,1]^d \cap s$, hence in $P$.  Therefore $P^\circ = s \cap (0,1]^d$, which clearly spans $s$.
\end{proof}

If $P$ is contained in a coordinate hyperplane, one should reformulate the problem with fewer variables.

Now for the main theorem on cubical magic.
The magic subspace $s \subseteq \bbR^d$ is given by \eqref{E:homoeqns}; the polytope is $P := [0,1]^d \cap s$.  For $t = 1,2,\ldots$ let 
\begin{enumerate}
\item[]
\begin{enumerate}
\item[$B^\circ(t) :=$] the number of integer points $x \in s$ with distinct entries that satisfy $0 < x_i < t$, 
\end{enumerate}
\end{enumerate}
and let 
\begin{enumerate}
\item[]
\begin{enumerate}
\item[$B(t) :=$] the number of pairs $(x,\sigma)$ consisting of an integer point $x \in s$ that satisfies $0 \leq x_i \leq t$ and a compatible $s$-realizable permutation $\sigma$ of $[d]$. 
\end{enumerate}
\end{enumerate}

\begin{thm}[Magic enumeration by bounds] \label{T:homomagic} 
Suppose $P := [0,1]^d \cap s$ does not lie within a coordinate hyperplane. 
Then $B^\circ$ and $B$ are quasipolynomials with leading term $(\vol P) t^{\dim s}$ and with constant term $B(0)$ equal to the number of permutations of $[d]$ that are realizable in $s$.  
Furthermore, $(-1)^{\dim s} B^\circ(-t) = B(t)$.
\end{thm} 

\begin{proof}
 Theorem \ref{T:magic} shows that $B = E_{\PH}$ and $B^\circ = 
E^\circ_\PoH$ are reciprocal quasipolynomials.  The constant term equals 
the number of permutations that are realizable in $P$.  The same ones 
are realizable in $s$, because a permutation realized by $x$ is also 
realized by $\alpha x$ for any positive real $\alpha$.  By choosing small 
enough $\alpha$ we can put $\alpha x$ into $P$ or $-P$.  If the latter, 
then $\frac12 \bj - \alpha x$ realizes the permutation and lies in $P$.
 \end{proof}

Most interesting is the case in which all the forms $f_i$ have weight zero. 
It is precisely then that $\cH$ and $P$ are transverse, as one can see by comparing the subspace $\bigcap\cH = \Span{\bj} \cap s$ with the definition of transversality; moreover, then $P$ spans $s$.

\begin{thm} \label{C:homomobius}
With $s$ defined by forms of weight zero, and assuming a magic labelling exists, we have
\begin{align*}
B^\circ(t) &= \sum_{u \in \cL(\PoH)} \mu(\0,u) E_{(0,1)^d \cap u}(t) \ , \\
B(t) &= \sum_{u \in \cL(\PoH)} |\mu(\0,u)| E_{[0,1]^d \cap u}(t) \ ,
\end{align*}
 where $\mu$ is the M\"obius function of $\cL(\PoH)$.
\end{thm}

\begin{proof}
Apply Theorem \ref{T:ehrhyp} in $s$.  A magic labelling exists if and only if $s$ does not lie in any hyperplane $x_j=x_k$.
\end{proof}

A flat $u \in \cL(\cH)$ has the form $v \cap s$ where $v$ is given by a 
series of equations of coordinates: $x_{i_1} = x_{i_2} = \cdots = x_{i_p}$, 
$x_{j_1} = x_{j_2} = \cdots = x_{j_q}$, etc.; that is, $v$ corresponds to a 
partition $\pi$ of $X$.  We can treat these equations as eliminating the 
variables $x_{i_2}, \ldots, x_{i_p}, x_{j_2}, \ldots, x_{j_q}, \ldots$ in 
favor of $x_{i_1}, x_{j_1}, \ldots$ .  With this substitution, $v = 
\bbR^{d'}$ for some $d' < d$ and $[0,1]^d \cap v = [0,1]^{d'}$.  Then 
$[0,1]^d \cap u$ is essentially $[0,1]^{d'} \cap s'$, where $s'$ is $s$ 
after identifying variables.  Similar remarks apply to $\cL(\PoH)$, since 
its members are the nonvoid intersections with $P^\circ$ of flats of 
$\cH$.

We describe some of the most interesting examples, concluding with the two best known.

\begin{exam}[Lines of constant length]\label{X:homolineswt}
In a covering clutter $(X,\cL)$, suppose every line has the same number of points.  Then the linear equations that express the existence of a magic sum take the form
$$
\sum_{j\in L_1} x_j = \sum_{j \in L_2} x_j \text{ for all } L_1, L_2 \in \cL \ .
$$
Theorems \ref{T:homomagic} and \ref{C:homomobius} both apply: $B^\circ(t)$ and $B(t)$ are reciprocal quasipolynomials in the upper bound $t-1$ or $t$, respectively, and so on.  Examples include magic, semimagic, and pandiagonal magic squares, affine and projective planes, $k$-nets, and magic hypercubes with or without diagonals of various kinds, as well as magic labelling of regular graphs.  One has to ask about the existence of a magic labelling.  There is no known general answer, but certainly there exist magic and semimagic squares of all orders $n \geq 3$ and pandiagonal magic squares of all orders $n \geq 4$, using the standard entries $\{1,2,\ldots,n^2\}$ if $n \not\equiv 2 \mod 4$ (see \cite[pp.\ 203--211]{Ball}).
\end{exam}

For a general covering clutter only Theorem \ref{T:homomagic} applies.  If it has magic labellings at all, the problems of existence and characterization of realizable permutations come to the fore.  
An example of nonexistence is the Fano plane.  Examples of existence are magic squares of size $n \geq 3$.  
Regarding characterization we propose a conjecture.  We may assume $X=[d]$.  
Given a covering clutter, a \emph{magic permutation} (with respect to the covering clutter) is a permutation of $[d]$ that is realizable by a \emph{positive} point $x \in P$.  Obviously, $x$ can be chosen to be rational if it exists at all.  In the cubical situation we are discussing here, all points in $P^\circ$ are positive, so magic permutations are identical to $P$-realizable and therefore to $s$-realizable permutations.  
A permutation $\sigma$ of $[d]$ defines a \emph{reverse dominance order} on the power set $\cP([d])$ by
\begin{enumerate}
\item[] $L \preccurlyeq_\sigma L'$ if, when $L$ and $L'$ are written in decreasing order according to $\sigma$, say $L = \{ \sigma j_1, \ldots, \sigma j_l \}$ where $j_1 > \cdots > j_l$ and $L' = \{ \sigma j'_1, \ldots, \sigma j'_{l'} \}$ where $j'_1 > \cdots > j'_{l'}$, then $l \leq l'$ and $j_1 \leq j'_1, \ldots, j_l \leq j'_l$.
\end{enumerate}
This is a partial order on $\cP([d])$.

\begin{conj}[Magic permutations] \label{Cj:magicperms}
A permutation $\sigma$ of $[d]$ is realizable by a positive point in the magic subspace $s$ of the covering clutter $([d],\cL)$ if and only if $\cL$ is an antichain in the reverse dominance order due to $\sigma$.
\end{conj} 

The archetype is magic squares.  A magic permutation there is a permutation of the cells of the square that is obtained from some magic labelling by arranging the cells in increasing order.  
We have verified the conjecture for $3\times3$ magic and semimagic squares (Examples \ref{X:magichomo3} and \ref{X:semimagichomo3}).

A permutation being a total order on $[d]$, one could generalize to total preorders (in which the antisymmetric law is not required); we propose the corresponding conjecture, replacing ``permutation'' by ``total preorder''.

The necessity for $\cL$ to be an antichain is obvious.  It is the same 
with the extension to linear forms.  We define the \emph{reverse dominance 
order due to $\sigma$} on the set $\left(\bbR^d\right)^*_{+}$ of positive 
linear forms on $\bbR^d$ by
\begin{enumerate}
\item[] $f \preccurlyeq_\sigma f'$ if, writing $f = \sum_{k=1}^d a_{k}x_{k}$ and $f' = \sum_{k=1}^d a'_{k}x_{k}$, then
$$
\sum_{k=j}^d a_{\sigma k} \leq \sum_{k=j}^d a'_{\sigma k} \text{ for every }j = 1,2,\ldots,d .
$$
\end{enumerate}

\begin{conj}[Magic permutations for linear forms] \label{Cj:magicformperms}
Let $f_1, \ldots, f_m$ be positive forms on $\bbR^d$ and $s$ the subspace on which they are all zero.  A permutation $\sigma$ of $[d]$ is realizable by a positive point in $s$ if and only if $\{ f_1, \ldots, f_m \}$ is an antichain in the reverse dominance order of forms.
\end{conj} 

\begin{exam}[Cubical symmetry] \label{X:homosym}
 A \emph{cubically symmetric} magic or semimagic square has the property 
that any two cells that lie opposite each other across the center have sum 
equal to $t$, the cell-value bound, and the center cell (if there is one) 
contains the value $\frac12 t$.  This definition is generalized from that 
of \emph{associated} square in \cite {Andrews} (\emph{symmetrical} square 
in \cite{Ball}), in which the entries are $1,2,\ldots,n^2$ and a 
symmetrical pair sums to $n^2+1$.  (Cubical symmetry contrasts with affine 
symmetry, which we shall treat shortly.)  The novel feature is the 
additional linear restraints besides the magic sum conditions: these are 
$y_{ij} + y_{n+1-i,n+1-j} = t$ for all $i,j \in [m]$.  Translated into the 
language of $\frac 1t$-fractional vectors $x = \frac 1t y \in \frac 1t 
\bbZ^{n^2}$, we require
\begin{equation} \label{E:homosym}
x_{ij} + x_{n+1-i,n+1-j} = 1.
\end{equation}
 The effect is to reduce the magic subspace $s$ to a smaller subspace $s'$, but Theorems \ref{T:magic} and \ref{C:homomobius} apply (with $s'$ replacing $s$) if the line size is constant, due to the next lemma.  In it we generalize symmetry to fairly arbitrary covering clutters on the point set $[n]^2$. 
The definition is the same: we take the magic sum conditions and the symmetry equations \eqref{E:homosym}.  
The magic subspace $s \subseteq \bbR^{n^2}$ is defined by equality of line sums; the \emph{cubically 
symmetric magic subspace}, $s'$, is the affine subspace of $s$ in which 
\eqref{E:homosym} is valid; then $P = [0,1]^{n^2} \cap s'$ and $\cH = \cH[K_{n^2}]^{s'}$.

\begin{lem} \label{L:homosym}
 If $\left( [n]^2, \cL \right)$ is a covering clutter of constant size such that $s'$ is not contained in a hyperplane $x_{ij} = x_{kl}$ of $\cH[K_{n^2}]$, and in particular if a cubically symmetric, strongly magic labelling exists, then $P$ and $\cH$ are transverse.
 \end{lem}

\begin{proof}
Since $\frac12 \bj \in P \cap \bigcap\cH$, every flat of $\cH$ intersects $P^\circ$.
\end{proof}

The hypothesis on $s'$ is satisfied in the case of cubically symmetric magic squares of side $n \geq 3$ because such squares are known to exist, using the values $1,2,\ldots,n^2$ if $n \not\equiv 2 \mod 4$; in fact, cubically symmetric pandiagonal squares exist if $n \geq 4$.  See \cite[pp.\ 204--211]{Ball}.
 \end{exam}

\begin{exam}[Magic squares of order 3: cubical count] \label{X:magichomo3}
 Define $M_\cc^\circ(t)$ for $t=1,2,3,\ldots$ to be the number of 
$3\times3$ magic squares in which each cell value is less than $t$. 
 From \cite{SLS}, the Ehrhart quasipolynomial is
 $$
M_\cc^\circ(t) = \begin{cases}
\frac{t^3-16t^2+76t-96}{6} = 
\frac{(t-2)(t-6)(t-8)}{6} &\text{if } t \equiv 0,2,6,8 \mod 12 , \\
 \\
\frac{t^3-16t^2+73t-58}{6} = 
\frac{(t-1)(t^2-15t+58)}{6} &\text{if } t \equiv 1 \mod 12 , \\
 \\
\frac{t^3-16t^2+73t-102}{6} = 
\frac{(t-3)(t^2-13t+34)}{6} &\text{if } t \equiv 3,11 \mod 12 , \\
 \\
\frac{t^3-16t^2+76t-112}{6} = 
\frac{(t-4)(t^2-12t+28)}{6} &\text{if } t \equiv 4,10 \mod 12 , \\
 \\
\frac{t^3-16t^2+73t-90}{6} = 
\frac{(t-2)(t-5)(t-9)}{6} &\text{if } t \equiv 5,9 \mod 12 , \\
 \\
\frac{t^3-16t^2+73t-70}{6} = 
\frac{(t-7)(t^2-9t+10)}{6} &\text{if } t \equiv 7 \mod 12 .
\end{cases}
$$

The constant term of $M_\cc(t) = (-1)^3 M^\circ_\cc(-t)$ is the number of magic permutations, which 
is 16.  These permutations are the rotations and reflections of the patterns
 $$
\text{(a)} \quad
\begin{tabular}{|c|c|c|}
\hline 4 & 9 & 2 \\ \hline 3 & 5 & 7 \\ \hline 8 & 1 & 6 \\ \hline
\end{tabular}
\hspace{1in}
\text{(b)} \quad
\begin{tabular}{|c|c|c|}
\hline 3 & 9 & 2 \\ \hline 4 & 5 & 6 \\ \hline 8 & 1 & 7 \\ \hline
\end{tabular}
$$
In these diagrams the numbers are not cell values but rather permutation positions: the largest value is in the cell marked 9, the next largest in that marked 8, and so on.  (To realize the permutations by magic squares, (a) can be left untouched but (b) needs numbers.)  The general form of a magic square is, up to the eight symmetries and an additive constant on each value,
$$
\begin{tabular}{|c|c|c|}
\hline 
\vstrut{-1.5ex}{4ex}\makebox[4em]{$-\beta$} & \makebox[4em]{$\alpha+\beta$} & \makebox[4em]{$-\alpha$} \\ 
 \hline
\vstrut{-1.5ex}{4ex}$-(\alpha-\beta)$ & $0$ & $\alpha-\beta$ \\ 
 \hline
\vstrut{-1.5ex}{4ex}$\alpha$ & $-(\alpha+\beta)$ & $\beta$ \\ 
 \hline
\end{tabular}
\vspace{10pt}
$$
where $\alpha > \beta > 0$ and $\alpha \neq 2\beta$.  If $\alpha > 2\beta$ 
we get the magic permutation (a); if $\alpha < 2\beta$ we get (b).  This proves 
there are just 16 magic permutations.
\end{exam}

\begin{exam}[Semimagic squares of order 3: cubical count]
 \label{X:semimagichomo3}
 Let $S_{\cc}^\circ(t)$, for $t>0$, be the number of semimagic squares of 
order 3 in which every entry is less than $t$.  \cite{SLS} has exact 
formulas. 
 The constant term $|S_{\cc}(0)| = 1296 = 6^4$ equals the number of $3\times3$ \emph{semimagic permutations}, that is, magic permutations for semimagic squares of order 3, in agreement with Conjecture \ref{Cj:magicperms}.
 We verified this by hand, finding all semimagic permutations of order 3, 
based on the fact that a normalized semimagic permutation is a linear 
extension of the partial ordering implied by the supernormalized form of a 
semimagic square developed in \cite{SLS}.
 \end{exam}

\subsection{Affine magic} \label{affmagic}

The affine approach counts magic squares, and magic labellings in general, by the magic sum.  In the general situation the magic subspace $s$ is defined by a rational, nonhomogeneous linear system
\begin{equation} \label{E:affeqns}
f_1(x) = f_2(x) = \cdots = f_m(x) = 1
\end{equation}
that we assume is consistent.  The magic polytope $P$ is the nonnegative part of $s$, that is,
$$
P := s \cap \Orth, 
$$
where $\Orth := \bbR^d_{\geq0}$, the nonnegative orthant, and the hyperplane arrangement is $\cH := \cH[K_d]^s$.  

Affine magic is quite similar to cubical magic, but there is something 
new: one has to worry about boundedness of $P$.  Obviously, $P$ is bounded 
if the defining linear forms $f_i$ in \eqref{E:affeqns} are positive and 
every variable appears in a form.  If the latter fails, we are simply in the wrong dimension, so we make the overall assumption in this section that \emph{every variable appears in at least one form}.  (A covering clutter satisfies this automatically.)  
As in the cubical treatment, one must make sure that $P$ affinely spans $s$ (or else change $s$ in the theorems to $\aff P$) and one has to be concerned about transversality of $P$ and $\cH$.

\begin{lem}\label{L:affspan}
If $P$ is not contained within a coordinate hyperplane, it spans $s$.
\end{lem}

\begin{proof}
As with Lemma \ref{L:homospan}.
\end{proof}

Let $s \subseteq \bbR^d$ be the solution space of \eqref{E:affeqns}, where 
the $f_i$ are rational linear forms, and let $P := s \cap \Orth$.  For $t 
= 1,2,\ldots$, let
\begin{enumerate}
\item[]
\begin{enumerate}
\item[$A^\circ(t) :=$] the number of integer points $x \in tP$ with distinct positive entries, 
\end{enumerate}
\end{enumerate}
and let 
\begin{enumerate}
\item[]
\begin{enumerate}
\item[$A(t) :=$] the number of pairs $(x,\sigma)$ consisting of a nonnegative integer point $x \in tP$ and a compatible $P^\circ$-realizable permutation $\sigma$ of $[d]$. 
\end{enumerate}
\end{enumerate}

When $P$ is bounded and spans $s$, Theorem \ref{T:magic} applies.

\begin{thm}[Magic enumeration by line sums] \label{T:affmagic} 
 Suppose that $P$ is bounded and does not lie within a coordinate 
hyperplane.  Then $A^\circ$ and $A$ are quasipolynomials with leading term 
$(\vol P) t^{\dim s}$ and with constant 
term $A(0)$ equal to the number of permutations of $[d]$ that are 
realizable in $P^\circ$.  Furthermore, $(-1)^{\dim s} A^\circ(-t) = A(t)$.
 \end{thm} 

\begin{proof}
A straightforward application of Theorem \ref{T:magic}.
\end{proof}

\begin{prob}\label{Pr:affs}
In Theorem \ref{T:affmagic}, can realizability in $P^\circ$ be replaced by the weaker property of realizability in $s$?
\end{prob}

We want to know when $P$ and $\cH$ are transverse so that the M\"obius-function formulas of Theorem \ref{T:ehrhyp} will apply, and also to allow a positive answer to the question of Problem \ref{Pr:affs}, for which the combination of transversality and nonemptiness of $P\cap\Span{\bj}$ is necessary and  sufficient.

\begin{lem} \label{L:afftrans}
$P$ and $\cH$ are transverse if $s$ is not contained in any hyperplane $x_j = x_k$ and all forms $f_i$ have equal positive weight.
\end{lem}

\begin{proof}
 We remarked that $P$ is bounded because every coordinate appears in a form.  
We need to verify that for no $u \in \cL(\cH)$ is $\eset \neq P\cap u \subseteq \partial P$.  But $P \cap u = P \cap v$ for some $v \in \cL(\cH[K_d])$, and $\partial P$ is contained in the union of the coordinate hyperplanes since they determine the facets of $P$.  If some $v \cap P$ lies in a coordinate hyperplane, then the same is true of $v = \Span{\bj}$; but $\Span{\bj} \cap P$ is in a coordinate hyperplane if and only if it is $\eset$ or $\{0\}$.  The latter is impossible with forms as in \eqref{E:affeqns}. 
 If the forms have equal weight $c > 0$, then $\Span{\bj} \cap P = \{ c\inv \bj \}$.  The converse is obvious.
\end{proof}

\begin{thm} \label{C:affmobius}
 With $s$ defined by forms of constant positive weight, and assuming a 
magic labelling exists, we have
\begin{align*}
 A^\circ(t) &= \sum_{u \in \cL(\PoH)} \mu(\0,u) E_{u \cap \Orth^\circ}(t) 
\ , \\
 A(t) &= \sum_{u \in \cL(\PoH)} |\mu(\0,u)| E_{u \cap \Orth}(t) \ ,
\end{align*}
 where $\mu$ is the M\"obius function of $\cL(\PoH)$. 
 \end{thm}

\begin{proof}
Transversality holds by Lemma \ref{L:afftrans}, since a magic labelling exists if and only if $s$ does not lie in any hyperplane $x_j=x_k$.  
Apply Theorem \ref{T:ehrhyp} in $s$.
\end{proof}

Some interesting examples are the affine versions of the cubical examples 
we already mentioned.  In many of them the affine intersection poset is 
combinatorially equivalent to the cubical intersection poset.  Let $P_\cc$ 
be the cubical polytope $s\cap(0,1)^d$; $P$ is still the affine polytope 
$s\cap\Orth$.

\begin{lem}\label{L:intersectionisom}
Suppose that $s$ is defined by positive forms such that each variable 
$x_j$ has some form in which its coefficient is at least $1$.  Assume also 
that a magic labelling exists.  Then 
the affine and cubical intersection posets are naturally isomorphic: 
$\cL(P_\cc^\circ,\cH[K_d]) \cong \cL(P^\circ,\cH[K_d])$ 
by $u\cap P_\cc^\circ \mapsto u\cap P^\circ$ for $u \in \cL(\cH[K_d])$ 
such that $u \cap P_\cc \neq \eset$.
 \end{lem}

\begin{proof}
 The assumption on the coefficients implies that $P = P_\cc \cap s$.
 \end{proof}

\begin{exam}[Lines of constant length, cf.\ Example \ref{X:homolineswt}] 
\label{X:afflineswt}
 The linear equations that express the existence of a magic sum $t$ take 
the form
$$
\sum_{j\in L} x_j = t \text{ for all } L \in \cL \ .
$$
 Theorems \ref{T:affmagic} and \ref{C:affmobius} and Lemma \ref{L:intersectionisom} all apply as long 
as a magic labelling exists; thus, $A^\circ(t)$ and $A(t)$ are reciprocal 
quasipolynomials in the magic sum $t$, etc.
\end{exam}

For general covering clutters Theorem \ref{T:affmagic} applies, showing that $A(t)$ and $A^\circ(t)$ are reciprocal quasipolynomials in $t$---provided that magic labellings exist at all.  With that assumption, Lemma \ref{L:intersectionisom} also applies.
The affine magic subspace, call it $s_1$, defined by $\sum_{i\in L} x_i = 1$ for all $L \in \cL$, is an affine subspace of the homogeneous magic subspace, call it $s_0$, and $s_0$ is the linear subspace generated by $s_1$.  This means that the permutations realizable in $s_0$ and $s_1$ are the same, so all our comments on magic permutations with respect to a covering clutter, in the context of cubical counting, apply as well to affine enumeration, except that we do not in general know that all $s$-realizable permutations are $P$-realizable.

\begin{exam}[Affine symmetry; cf.\ Example \ref{X:homosym}] \label{X:affsym}
An \emph{affinely symmetric} magic or semimagic square has the property that the average value of any two cells that lie opposite each other across the center equals the average cell value, $t/n$, and the center cell (if there is one) contains the value $t/n$.  (Of course, one cannot expect such squares to exist unless $t \equiv 0 \pmod{n}$.)  This definition is another generalization of that of associated square.  The additional linear restraints are 
$
y_{ij} + y_{n+1-i,n+1-j} = 2t/n \text{ for all } i,j \in [m].  
$
Translated into the language of $\frac 1t$-fractional vectors $x = \frac 1t y \in \frac 1t \bbZ^{n^2}$, we require
\begin{equation} \label{E:affsym}
x_{ij} + x_{n+1-i,n+1-j} = \frac 2 n \ .
\end{equation}
The effect is to reduce the magic subspace $s$ to a smaller subspace $s'$.  
With the extra hypothesis that all lines have $n$ points, Theorems \ref{T:magic} and \ref{C:affmobius} apply, with $s'$ replacing $s$, by the following lemma, in which we generalize symmetry to covering clutters on $X = [n]^2$, with magic sum conditions and the symmetry equations \eqref{E:affsym}.  
The magic subspace $s \subseteq \bbR^{n^2}$ is defined by equality of line sums; the \emph{affinely symmetric magic subspace}, $s'$, is the affine subspace of $s$ in which \eqref{E:homosym} is valid; then $P = \bbR_{\geq 0}^{n^2} \cap s'$ and $\cH = \cH[K_{n^2}]^{s'}$.

\begin{lem} \label{L:affsym}
If $\left( [n]^2, \cL \right)$ is a covering clutter of size $n$ such that $s'$ is not contained in a hyperplane $x_{ij} = x_{kl}$ of $\cH[K_{n^2}]$, and in particular if an affinely symmetric magic labelling exists, then $P$ and $\cH$ are transverse.
\end{lem}

\begin{proof}
Since $\frac 1n \bj \in P \cap \bigcap\cH$, every flat of $\cH$ intersects $P^\circ$.
\end{proof}
\end{exam}

 \begin{exam}[Magic squares of order 3: affine count]
\label{X:magicaff3}
 Let $M_\aa^\circ(t)$, for $t>0$, be the number of magic $3\times3$ 
squares with magic sum $t$.  In \cite{SLS} we find that 
 $$
M_\aa^\circ(t) = \begin{cases}
\frac{2t^2-32t+144}{9} = 
\frac{2}{9}(t^2-16t+72)  &\text{if } t \equiv 0 \mod 18 , \\
 \\
\frac{2t^2-32t+78}{9} = 
\frac{2}{9}(t-3)(t-13) &\text{if } t \equiv 3 \mod 18 , \\
 \\
\frac{2t^2-32t+120}{9} = 
\frac{2}{9}(t-6)(t-10) &\text{if } t \equiv 6 \mod 18 , \\
 \\
\frac{2t^2-32t+126}{9} = 
\frac{2}{9}(t-7)(t-9) &\text{if } t \equiv 9 \mod 18 , \\
 \\
\frac{2t^2-32t+96}{9} = 
\frac{2}{9}(t-4)(t-12) &\text{if } t \equiv 12 \mod 18 , \\
 \\
\frac{2t^2-32t+102}{9} = 
\frac{2}{9}(t^2-16t+51)  &\text{if } t \equiv 15 \mod 18 , \\
 \\
0  &\text{if } t \not\equiv 0 \mod 3 .
\end{cases}
 $$

Xin \cite{Xin} has another way to find the generating function of $M_\aa^\circ(t)$ (with the minor difference that he allows zero entries), by applying MacMahon's partition calculus.

The constant term $M_\aa(0) = 16$ is the number of magic permutations; 
this is the same number as with cubically counted magic squares, 
$M_\cc(0)$ in Example \ref{X:magichomo3}.
 \end{exam}

\begin{exam}[Semimagic squares of order 3: affine count] \label{X:semimagicaff3}
 Complete results for this example are in \cite{SLS}.
 \end{exam}

\section{Latin squares join in magically} \label{latin}

The general picture that encompasses latin squares is that of a covering clutter $(X,\cL)$, as in Section \ref{magic}, with an integer labelling $x: X \to \bbZ$ subject to the requirement that 
$$
x(e) \neq x(f) \qquad \text{ if $e$ and $f$ lie in a line.}
$$
This is a \emph{latin labelling} of $(X,\cL)$.  The graph of forbidden equalities is therefore
$$
\Gamma_{\cL} := \bigcup_{L \in \cL} K_L,
$$
$K_L$ being the complete graph with $L$ as node set.  Every graph is equal to $\Gamma_\cL$ for some choice of covering clutter.  An orientation of $\Gamma_{\cL}$ is \emph{acyclic} if it has no directed cycles.  An orientation and a node $c$-coloring $x : V\to [c]$ are \emph{compatible} if $x_j\geq x_i$ whenever there is a $\Gamma_{\cL}$-edge oriented from $i$ to $j$ \cite{AOG}.

A crucial decision is how to restrict the symbols of the latin labelling.  
One may simply specify the number of symbols allowed, say $k$, and an 
arbitrary symbol set, let us say $[k]$.  Then the number of latin 
labellings equals the chromatic polynomial $\chi_{\Gamma_\cL}(k)$, so we would have merely an application of graph coloring. 
 In this article, taking a leaf from the guidebook of magic squares, we focus on approaches which give rise to \emph{magilatin squares}: we impose a summation condition on the lines, which may be either homogeneous: 
\begin{enumerate} 
\item[(i)] Set all line sums equal to each other.
\end{enumerate} 
or affine: 
\begin{enumerate} 
\item[(ii)] Set all line sums equal to $t$.
\end{enumerate} 
 This tactic is not so strange as it may appear.  
A partial magilatin square with homogeneous line-sum requirements and symbols restricted to the interval $[1,n]$ is just a latin square with symbol set $[n]$.
So is a partial latin square with affine line-sum restrictions, line sum $t = \binom{n+1}{2}$, and positive symbols.  
We are assuming that the affine constraint (ii) is supplemented by a positivity assumption and that the homogeneous constraint (i) is supplemented by the requirement that the symbols be drawn from the set $[k]$ for some $k$.
To handle partial latin rectangles requires a generalization of (i) to multiple covering clutters; see Section \ref{homolatin}.  
In every case the hyperplane arrangement is $\cH := \cH[\Gamma_\cL]^s$, where $s$ is the subspace of $\bbR^X$ determined by the appropriate line sum conditions.

\subsection{Cubical latinity} \label{homolatin}

The cubical approach to latin labellings counts them by the largest allowed value; if $t$ is the parameter, the labellings counted are those with $0 < x_i < t$.  
We assume a \emph{multiple covering clutter}, $(X;\cL_1,\ldots,\cL_k)$, that is, there are $k\geq1$ covering clutters $\cL_1,\ldots,\cL_k$; in our examples $k$ will be $1$ for squares and $2$ for rectangles.  The labelling must be such that the lines in each class have equal sums; but the sums in different classes are independent of each other.  For a magilatin labelling the numbers in each line must be distinct.
To count the labellings we take the subspace $s$ of $\bbR^d$ in which all line sums within each covering clutter $\cL_i$ are equal.  The polytope is $P := [0,1]^d \cap s$.  This is as in Section \ref{homomagic}, so Lemma \ref{L:homospan} applies to assure that $P$ spans $s$.  For $t = 1,2,\ldots$, let 
\begin{enumerate}
\item[]
\begin{enumerate}
\item[$L_\cc^\circ(t) :=$] the number of latin labellings $x$ with equal line sums within each covering clutter and with entries that satisfy $0 < x_i < t$, 
\end{enumerate}
\end{enumerate}
and let 
\begin{enumerate}
\item[]
\begin{enumerate}
\item[$L_\cc(t) :=$] the number of pairs consisting of a latin labelling $x$, with equal line sums within each covering clutter and with $0 \leq x_i \leq t$, and a compatible acyclic orientation of $\Gamma_\cL$. 
\end{enumerate}
\end{enumerate}

\begin{thm}[Magilatin enumeration by upper bound] \label{T:homolatin} 
Suppose $P := [0,1]^d \cap s$ does not lie within a coordinate hyperplane.  
Then $L_\cc^\circ$ and $L_\cc$ are quasipolynomials with leading term $(\vol P) t^{\dim s}$ and with constant term $L_\cc(0)$ equal to the number of acyclic orientations of $\Gamma_\cL$ that are realizable in $s$.  
Furthermore, $(-1)^{\dim s} L_\cc^\circ(-t) = L_\cc(t)$.
\end{thm} 

\begin{proof}
Theorem \ref{T:quasi} shows that $L_\cc = E_{\PH}$ and $L_\cc^\circ = E^\circ_\PoH$ are reciprocal quasipolynomials.  The remainder of the proof is as in Theorem \ref{T:homomagic}.
\end{proof}

\begin{thm} \label{C:homolatinmobius}
 If all lines have the same size and a latin labelling exists, then we 
have all the formulas of Theorem \ref{C:homomobius} with $L_\cc$ in 
place of $B$.
 \end{thm}

\begin{proof}
Apply Theorem \ref{T:ehrhyp}.
\end{proof}

\begin{exam}[Small latin shapes, counted cubically] \label{X:latinhomo}
We calculated two very small magilatin examples: a square and a rectangle.

In the $2\times2$ magilatin square it must be that $x_{11}=x_{22} \neq x_{12}=x_{21}$.  With cubical constraints, then, 
$$
L_\cc^\circ(t) = (t-1)(t-2) .
$$
For comparison, the number of squares without the magilatin distinctness requirement is $t^2$.

The $2\times3$ latin rectangle is much more complicated.  Calculation with Maple shows that the denominator of $P$ is 6 and that of $(\PH)$ is 12.  
This permits us to calculate the quasipolynomial (by actual count of rectangles up to $t=60$).  Its period turns out to be 4.
$$
L_\cc^\circ(t) = \begin{cases}
\frac{t^3-12t^2+41t-30}{4} = \frac{(t-1)(t-5)(t-6)}{4} &\text{if } t \equiv 1 \mod 4, \\
\frac{t^3-12t^2+41t-42}{4} = \frac{(t-2)(t-3)(t-7)}{4} &\text{if } t \equiv 3 \mod 4, \\
\frac{t^3-12t^2+44t-48}{4} = \frac{(t-2)(t-4)(t-6)}{4} &\text{if $t$ is even}.
\end{cases}
$$
For comparison, the number of rectangles without the distinctness requirement is
$$
\begin{cases}
\frac{t^3-3t^2+3t-1}{4} = \frac{(t-1)^3}{4} &\text{if $t$ is odd}, \\
\frac{t^3-3t^2+6t-4}{4} = \frac{(t-1)(t^2-2t+4)}{4} &\text{if $t$ is even}.
\end{cases}
$$

\begin{prob}\label{Pr:latinpolyform}
 We do not know why the coefficients alternate in sign, what causes the smallness of the differences among the constituents, nor what makes the even constituents have a smaller period than the odd constituents.
 \end{prob}
\end{exam}

\begin{exam}[Magilatin $3\times3$ squares, counted cubically]
\label{X:latinhomo3}
See \cite{SLS} for the complete solution.
\end{exam}

\subsection{Affine latinity} \label{affinelatin}

For affine counting of latin labellings we take a covering clutter $(X,\cL)$.  $s$ is the subspace in which all line sums equal 1, and $P = s \cap \Orth$, as in affine magic (Section \ref{affmagic}).  
As there, we can apply Lemma \ref{L:affspan} to conclude that $P$ spans $s$ in all interesting cases.  
One advantage over magic is that, because we consider only line sums and not general linear forms, $P$ is certain to be bounded.
For $t = 1,2,\ldots$, let 
\begin{enumerate}
\item[]
\begin{enumerate}
\item[$L_\aa^\circ(t) :=$] the number of latin labellings of $(X,\cL)$ with positive entries and all line sums equal to $t$, 
\end{enumerate}
\end{enumerate}
and let 
\begin{enumerate}
\item[]
\begin{enumerate}
\item[$L_\aa(t) :=$] the number of pairs consisting of a nonnegative latin labelling with all line sums equal to $t$ and a compatible acyclic orientation of $\Gamma_\cL$ that are realizable in $P^\circ$. 
\end{enumerate}
\end{enumerate}

\begin{lem} \label{L:afftranslatin}
$P$ and $\cH$ are transverse if a positive latin labelling exists.
\end{lem}

\begin{proof}
As with Lemma \ref{L:afftrans}.  The existence of a positive latin labelling implies that $s$ is not contained in any hyperplane $x_j = x_k$ for which $e_i$ and $e_j$ are collinear and, for every $v \in \cL(\cH[\Gamma_\cL])$, $v \cap P$ is not contained in a coordinate hyperplane.
\end{proof}

\begin{thm}[Magilatin enumeration by magic line sum] \label{T:affmagiclatin} 
Suppose that a positive latin labelling exists for some $t$.  Then $L_\aa^\circ$ and $L_\aa$ are quasipolynomials with leading term $(\vol P) t^{\dim s}$ and with constant term $L_\aa(0)$ equal to the number of acyclic orientations of $\Gamma_\cL$ that are realizable in $P$.  
Furthermore, $(-1)^{\dim s} L_\aa^\circ(-t) = L_\aa(t)$.
 \end{thm} 

\begin{proof}
By Theorem \ref{T:magic}, adapted to the latin nonequalities, together with Lemma \ref{L:afftranslatin} to ensure transversality so that realizability in $P$ is equivalent to realizability in $P^\circ$.
\end{proof}

If every line has the same size, then the acyclic orientations that are realizable in $P^\circ$ are the same as those realizable in $s$, because then $\Span{\bj}$ intersects $P^\circ$; see the discussion following Theorem \ref{T:affmagic}.

\begin{thm} \label{C:affmobiuslatin}
 Assuming a positive latin labelling exists and all lines have the same size, we have the formulas of Theorem \ref{C:affmobius} with $L_\aa$ in place of $A$.
 \end{thm}

\begin{proof}
Transversality holds by Lemma \ref{L:afftranslatin}.  
Apply Theorem \ref{T:ehrhyp} in $s$.
\end{proof}

\begin{exam}[Small magilatin squares, counted affinely]\label{X:latinaff2}
It is easy to see that in a $2\times2$ magilatin square with affine constraints,
 $$
 L_\aa^\circ(t) = \begin{cases}
t-1 &\text{if $t$ is odd}, \\
t-2 &\text{if $t$ is even}.
\end{cases}
 $$
 The period of $2$ equals the denominator of $(\PH)$.  
The number of positive squares without the magilatin distinctness requirement is $t-1$.
\end{exam}

\begin{exam}[Magilatin $3\times3$ squares, counted affinely] \label{X:latinaff3}
For the number of these squares see \cite{SLS}.
\end{exam}

\section{Generalized exclusions} \label{extensions}

We concentrated our treatment on exclusions of magic and latin type: that is, where all values, or all collinear values, are unequal.  
Other kinds of exclusion are possible.  We want to mention the very natural \emph{complementation} restrictions.  In cubical enumeration we call $x_i$ and $x_j$ \emph{complementary} if $x_i + x_j = t$.  If we forbid certain pairs of values to be complementary, we pass from graphs to signed graphs, since the rule $x_i + x_j \neq t$ corresponds to a negative edge $-ij$.  An inequality $x_i \neq x_j$ corresponds to a positive edge, $+ij$; thus ordinary edges are positive.  The exact application of signed graphs involves translation and halving of the centrally symmetric polytope $[-1,1]^d$ to $[0,1]^d$, along with corresponding translation of the signed-graphic hyperplane arrangement, as explained in \cite[Section 5]{IOP}.  
Because signed-graphic hyperplanes, as translated, give half-integral vertices, we expect a counting quasipolynomial with nonequalities and noncomplementarities given by a signed graph to have twice the period of a counting quasipolynomial pertaining to a similar unsigned graph of inequalities.  This is necessarily vague; we intend only a suggestion for research that we invite readers to explore.


\end{document}